\documentclass[a4paper,10pt]{amsart}

\usepackage{amsmath,amssymb,amsthm,amscd,latexsym,stmaryrd}
\usepackage[all]{xy}

\newcommand{\N}{{\ensuremath{\mathbb N}}}
\newcommand{\Z}{{\ensuremath{\mathbb Z}}}

\newcommand{\K}{{\ensuremath{\mathbb K}}}

\newcommand{\Ocal}{\ensuremath{\mathcal{O}}}
\newcommand{\Pcal}{\ensuremath{\mathcal{P}}}
\newcommand{\Qcal}{\ensuremath{\mathcal{Q}}}
\newcommand{\Lcal}{\ensuremath{\mathcal{L}}}
\newcommand{\Ecal}{\ensuremath{\mathcal{E}}}

\newcommand{\s}{\ensuremath{\mathbb{S}}}
\newcommand{\smodule}{\ensuremath{\mathbb{S}}\text{-module}}
\newcommand{\smodules}{\ensuremath{\mathbb{S}}\text{-modules}}

\newcommand{\bbone}{1\hspace{-2.6pt}\mbox{\normalfont{l}}}

\newcommand{\Com}{\mathcal{C}om}
\newcommand{\Lie}{\mathcal{L}ie}
\newcommand{\As}{\mathcal{A}s}
\newcommand{\Perm}{\mathcal{P}erm}
\newcommand{\Prelie}{\mathcal{P}re\mathcal{L}ie}
\newcommand{\Dias}{\mathcal{D}ias}
\newcommand{\Dend}{\mathcal{D}end}
\newcommand{\Tridend}{\mathcal{T}ri\mathcal{D}end}
\newcommand{\Trias}{\mathcal{T}rias}
\newcommand{\Comtrias}{\mathcal{C}om\mathcal{T}rias}
\newcommand{\Postlie}{\mathcal{P}ost\mathcal{L}ie}

\newcommand{\Comtwo}{\tc{\Com}}
\newcommand{\Lietwo}{\lc{\Lie}}
\newcommand{\Astwotc}{\tc{\As}} 
\newcommand{\Astwolc}{\lc{\As}} 
\newcommand{\Permtwo}{\tc{\Perm}}
\newcommand{\Prelietwo}{\lc{\Prelie}}
\newcommand{\Diastwo}{\tc{\Dias}}
\newcommand{\Dendtwo}{\lc{\Dend}}
\newcommand{\Triastwo}{\tc{\Trias}}
\newcommand{\Tridendtwo}{\lc{\Tridend}}
\newcommand{\Comtriastwo}{\tc{\Comtrias}}
\newcommand{\Postlietwo}{\lc{\Postlie}}

\renewcommand{\wp}[1]{\Pi_{#1}^{\mbox{\tiny w}}}

\newcommand{\atom}[4]{\alpha^{#1,#2}_{#3,#4}}
\newcommand{\atomb}[3]{\alpha^{#1}_{#2,#3}}
\newcommand{\maximal}[2]{\mu_{#1,#2}}
\newcommand{\maximalb}[1]{\mu_{#1}}

\newcommand{\aao}[1]{\widetilde{#1}} 

\newcommand{\Free}{\mathcal{F}}
\newcommand{\QO}[2]{\Free(#1)/(#2)}

\newcommand{\lc}[1]{\mbox{$#1^2$}}
\newcommand{\tc}[1]{\mbox{$^2#1$}}
\newcommand{\subtc}[1]{\mbox{$\scriptstyle{^{\scriptscriptstyle{2}}#1}$}}

\newcommand{\cz}{\scriptscriptstyle{\lor}}

\DeclareMathOperator{\sgn}{sgn}
\DeclareMathOperator{\weight}{w}
\DeclareMathOperator{\nblocks}{n}
\DeclareMathOperator{\op}{op}
\DeclareMathOperator{\Hom}{Hom}

\numberwithin{equation}{section}

\newtheorem{prop}[equation]{Proposition}
\newtheorem{thm}[equation]{Theorem}
\newtheorem{cor}[equation]{Corollary}

\theoremstyle{definition}
\newtheorem{rem}[equation]{Remark}
\newtheorem{defn}[equation]{Definition}
\newtheorem{ex}[equation]{Example}

\newtheorem{defnalpha}{Definition}

\newcommand{\iso}{\cong}

\newcommand{\bhs}{\hspace{10pt}}
\newcommand{\shs}{\hspace{5pt}}

\newcommand{\B}[3][{}]{\ensuremath{
 \xygraph{
!{<0pt,0pt>;<4pt,0pt>:<0pt,-4pt>::}
!{(0,3)}*{\scriptscriptstyle #1}
!{(0,2)}="a"
!{(0,0)}*{\scriptscriptstyle{\bullet}}="b"
!{(-2,-2)}="c"
!{(2,-2)}="d"
!{(-2,-3)}*{\scriptscriptstyle #2}
!{(2,-3)}*{\scriptscriptstyle #3}
"a"-"b"
"b"-"c"
"b"-"d"
}
}}

\newcommand{\W}[3][{}]{\ensuremath{
 \xygraph{
!{<0pt,0pt>;<4pt,0pt>:<0pt,-4pt>::}
!{(0,3)}*{\scriptscriptstyle #1}
!{(0,2)}="a"
!{(0,0)}*{\scriptscriptstyle{\circ}}="b"
!{(-2,-2)}="c"
!{(2,-2)}="d"
!{(-2,-3)}*{\scriptscriptstyle #2}
!{(2,-3)}*{\scriptscriptstyle #3}
"a"-"b"
"b"-"c"
"b"-"d"
}
}}

\newcommand{\Bczech}[3][{}]{\ensuremath{
 \xygraph{
!{<0pt,0pt>;<4pt,0pt>:<0pt,-4pt>::}
!{(0,3)}*{\scriptscriptstyle #1}
!{(0,2)}="a"
!{(0,0)}*{\scriptscriptstyle{\bullet}}="b"
!{(1.5,0)}*{\scriptscriptstyle{^{\cz}}}
!{(-2,-2)}="c"
!{(2,-2)}="d"
!{(-2,-3)}*{\scriptscriptstyle #2}
!{(2,-3)}*{\scriptscriptstyle #3}
"a"-"b"
"b"-"c"
"b"-"d"
}
}}

\newcommand{\Wczech}[3][{}]{\ensuremath{
 \xygraph{
!{<0pt,0pt>;<4pt,0pt>:<0pt,-4pt>::}
!{(0,3)}*{\scriptscriptstyle #1}
!{(0,2)}="a"
!{(0,0)}*{\scriptscriptstyle{\circ}}="b"
!{(1.5,0)}*{\scriptscriptstyle{^{\cz}}}
!{(-2,-2)}="c"
!{(2,-2)}="d"
!{(-2,-3)}*{\scriptscriptstyle #2}
!{(2,-3)}*{\scriptscriptstyle #3}
"a"-"b"
"b"-"c"
"b"-"d"
}
}}

\newcommand{\BB}[4][{}]{\ensuremath{
 \xygraph{
!{<0pt,0pt>;<4pt,0pt>:<0pt,-4pt>::}
!{(1,4)}*{\scriptscriptstyle #1}
!{(1,3)}="a"
!{(1,1)}*{\scriptscriptstyle{\bullet}}="b"
!{(-1,-1)}*{\scriptscriptstyle{\bullet}}="c"
!{(3,-1)}="d"
!{(3,-2)}*{\scriptscriptstyle #4}
!{(-3,-3)}="e"
!{(-3,-4)}*{\scriptscriptstyle #2}
!{(1,-3)}="f"
!{(1,-4)}*{\scriptscriptstyle #3}
"a"-"b"
"b"-"c"
"b"-"d"
"c"-"e"
"c"-"f"
}
}}

\newcommand{\WW}[4][{}]{\ensuremath{
 \xygraph{
!{<0pt,0pt>;<4pt,0pt>:<0pt,-4pt>::}
!{(1,4)}*{\scriptscriptstyle #1}
!{(1,3)}="a"
!{(1,1)}*{\scriptscriptstyle{\circ}}="b"
!{(-1,-1)}*{\scriptscriptstyle{\circ}}="c"
!{(3,-1)}="d"
!{(3,-2)}*{\scriptscriptstyle #4}
!{(-3,-3)}="e"
!{(-3,-4)}*{\scriptscriptstyle #2}
!{(1,-3)}="f"
!{(1,-4)}*{\scriptscriptstyle #3}
"a"-"b"
"b"-"c"
"b"-"d"
"c"-"e"
"c"-"f"
}
}}

\newcommand{\BW}[4][{}]{\ensuremath{
 \xygraph{
!{<0pt,0pt>;<4pt,0pt>:<0pt,-4pt>::}
!{(1,4)}*{\scriptscriptstyle #1}
!{(1,3)}="a"
!{(1,1)}*{\scriptscriptstyle{\circ}}="b"
!{(-1,-1)}*{\scriptscriptstyle{\bullet}}="c"
!{(3,-1)}="d"
!{(3,-2)}*{\scriptscriptstyle #4}
!{(-3,-3)}="e"
!{(-3,-4)}*{\scriptscriptstyle #2}
!{(1,-3)}="f"
!{(1,-4)}*{\scriptscriptstyle #3}
"a"-"b"
"b"-"c"
"b"-"d"
"c"-"e"
"c"-"f"
}
}}

\newcommand{\WB}[4][{}]{\ensuremath{
 \xygraph{
!{<0pt,0pt>;<4pt,0pt>:<0pt,-4pt>::}
!{(1,4)}*{\scriptscriptstyle #1}
!{(1,3)}="a"
!{(1,1)}*{\scriptscriptstyle{\bullet}}="b"
!{(-1,-1)}*{\scriptscriptstyle{\circ}}="c"
!{(3,-1)}="d"
!{(3,-2)}*{\scriptscriptstyle #4}
!{(-3,-3)}="e"
!{(-3,-4)}*{\scriptscriptstyle #2}
!{(1,-3)}="f"
!{(1,-4)}*{\scriptscriptstyle #3}
"a"-"b"
"b"-"c"
"b"-"d"
"c"-"e"
"c"-"f"
}
}}

\newcommand{\WWB}[5][{}]{\ensuremath{
 \xygraph{
!{<0pt,0pt>;<4pt,0pt>:<0pt,-4pt>::}
!{(2,5)}*{\scriptscriptstyle #1}
!{(2,4)}="a"
!{(2,2)}*{\scriptscriptstyle{\bullet}}="b"
!{(0,0)}*{\scriptscriptstyle{\circ}}="c"
!{(4,-0)}="d"
!{(4,-1)}*{\scriptscriptstyle #5}
!{(-2,-2)}*{\scriptscriptstyle{\circ}}="e"
!{(2,-2)}="f"
!{(2,-3)}*{\scriptscriptstyle #4}
!{(-4,-4)}="g"
!{(0,-4)}="h"
!{(-4,-5)}*{\scriptscriptstyle #2}
!{(0,-5)}*{\scriptscriptstyle #3}
"a"-"b"
"b"-"c"
"b"-"d"
"c"-"e"
"c"-"f"
"e"-"g"
"e"-"h"
}
}}

\newcommand{\WBB}[5][{}]{\ensuremath{
 \xygraph{
!{<0pt,0pt>;<4pt,0pt>:<0pt,-4pt>::}
!{(2,5)}*{\scriptscriptstyle #1}
!{(2,4)}="a"
!{(2,2)}*{\scriptscriptstyle{\bullet}}="b"
!{(0,0)}*{\scriptscriptstyle{\bullet}}="c"
!{(4,-0)}="d"
!{(4,-1)}*{\scriptscriptstyle #5}
!{(-2,-2)}*{\scriptscriptstyle{\circ}}="e"
!{(2,-2)}="f"
!{(2,-3)}*{\scriptscriptstyle #4}
!{(-4,-4)}="g"
!{(0,-4)}="h"
!{(-4,-5)}*{\scriptscriptstyle #2}
!{(0,-5)}*{\scriptscriptstyle #3}
"a"-"b"
"b"-"c"
"b"-"d"
"c"-"e"
"c"-"f"
"e"-"g"
"e"-"h"
}
}}

\newcommand{\WBWvar}[5][{}]{\ensuremath{
 \xygraph{
!{<0pt,0pt>;<4pt,0pt>:<0pt,-4pt>::}
!{(-1.5,5)}*{\scriptscriptstyle #1}
!{(-1.5,4)}="a"
!{(-1.5,2)}*{\scriptscriptstyle{\circ}}="b"
!{(-3.5,0)}="c"
!{(1.5,0)}*{\scriptscriptstyle{\bullet}}="d"
!{(-3.5,-2)}="e"
!{(-0.5,-2)}*{\scriptscriptstyle{\circ}}="f"
!{(3.5,-2)}="g"
!{(-3.5,-3)}*{\scriptscriptstyle #2}
!{(3.5,-3)}*{\scriptscriptstyle #5}
!{(-2.5,-4)}="h"
!{(1.5,-4)}="i"
!{(-2.5,-5)}*{\scriptscriptstyle #3}
!{(1.5,-5)}*{\scriptscriptstyle #4}
"a"-"b"
"b"-"c"
"b"-"d"
"c"-"e"
"d"-"f"
"d"-"g"
"f"-"h"
"f"-"i"
}
}}

\newcommand{\Bsmall}{\ensuremath{ 
 \xygraph{
!{<0pt,0pt>;<2.6pt,0pt>:<0pt,-2.6pt>::}
!{(0,2)}="a"
!{(0,0)}*{\scriptscriptstyle{\bullet}}="b"
!{(-2,-2)}="c"
!{(2,-2)}="d"
"a"-"b"
"b"-"c"
"b"-"d"
}
}}

\newcommand{\Wsmall}{\ensuremath{ 
 \xygraph{
!{<0pt,0pt>;<2.6pt,0pt>:<0pt,-2.6pt>::}
!{(0,2)}="a"
!{(0,0)}*{\scriptscriptstyle{\circ}}="b"
!{(-2,-2)}="c"
!{(2,-2)}="d"
"a"-"b"
"b"-"c"
"b"-"d"
}
}}

\newcommand{\Id}[1]{
\xygraph{
!{<0pt,0pt>;<4pt,0pt>:<0pt,-4pt>::}
!{(0,2)}="a"
!{(0,-2)}="b"
!{(0,-3)}*{\scriptscriptstyle #1}
"a"-"b"
}}

\newcommand{\WBcorolla}{
 \xygraph{
!{<0pt,0pt>;<4pt,0pt>:<0pt,-4pt>::}
!{(4,6)}="a"
!{(4,4)}*{\scriptscriptstyle{\bullet}}="b"
!{(3.25,3.25)}*{.}="bda"
!{(2.5,2.5)}*{.}="bdb"
!{(1.75,1.75)}*{.}="bdc"
!{(6,2)}="c"
!{(1,1)}*{\scriptscriptstyle{\bullet}}="d"
!{(-1,-1)}*{\scriptscriptstyle{\circ}}="e"
!{(-1.75,-1.75)}*{.}="eha"
!{(-2.5,-2.5)}*{.}="ehb"
!{(-3.25,-3.25)}*{.}="ehc"
!{(3,-1)}="f"
!{(1,-3)}="g"
!{(-4,-4)}*{\scriptscriptstyle{\circ}}="h"
!{(-6,-6)}="i"
!{(-2,-6)}="j"
!{(6,1)}*{\scriptscriptstyle{n}}="C" 
!{(4,-2)}*{\scriptscriptstyle{i+2}}="F" 
!{(1,-4)}*{\scriptscriptstyle{i+1}}="G" 
!{(-2,-7)}*{\scriptscriptstyle{2}}="J" 
!{(-6,-7)}*{\scriptscriptstyle{1}}="I" 
"a"-"b"
"b"-"c"
"d"-"f"
"d"-"e"
"e"-"g"
"h"-"i"
"h"-"j"
}
}

\newcommand{\XX}[6][{}]{\ensuremath{
 \xygraph{
!{<0pt,0pt>;<4pt,0pt>:<0pt,-4pt>::}
!{(1,4)}*{\scriptscriptstyle #1}
!{(1,3)}="a"
!{(1,1)}*{\scriptscriptstyle{\circ}}="b"
!{(2,1)}*{\scriptscriptstyle #6}
!{(-1,-1)}*{\scriptscriptstyle{\circ}}="c"
!{(0,-1)}*{\scriptscriptstyle #5}
!{(3,-1)}="d"
!{(3,-2)}*{\scriptscriptstyle #4}
!{(-3,-3)}="e"
!{(-3,-4)}*{\scriptscriptstyle #2}
!{(1,-3)}="f"
!{(1,-4)}*{\scriptscriptstyle #3}
"a"-"b"
"b"-"c"
"b"-"d"
"c"-"e"
"c"-"f"
}
}}

\newcommand{\YY}[6][{}]{\ensuremath{
 \xygraph{
!{<0pt,0pt>;<4pt,0pt>:<0pt,-4pt>::}
!{(1,4)}*{\scriptscriptstyle #1}
!{(1,3)}="a"
!{(1,1)}="b"
!{(2,1)}*{\scriptscriptstyle #6}
!{(-1,-1)}="c"
!{(0,-1)}*{\scriptscriptstyle #5}
!{(3,-1)}="d"
!{(3,-2)}*{\scriptscriptstyle #4}
!{(-3,-3)}="e"
!{(-3,-4)}*{\scriptscriptstyle #2}
!{(1,-3)}="f"
!{(1,-4)}*{\scriptscriptstyle #3}
"a"-"b"
"b"-"c"
"b"-"d"
"c"-"e"
"c"-"f"
}
}}

\newcommand{\XZ}[6][{}]{\ensuremath{
 \xygraph{
!{<0pt,0pt>;<4pt,0pt>:<0pt,-4pt>::}
!{(1,4)}*{\scriptscriptstyle #1}
!{(1,3)}="a"
!{(1,1)}*{\scriptscriptstyle{\bullet}}="b"
!{(2,1)}*{\scriptscriptstyle #6}
!{(-1,-1)}*{\scriptscriptstyle{\circ}}="c"
!{(0,-1)}*{\scriptscriptstyle #5}
!{(3,-1)}="d"
!{(3,-2)}*{\scriptscriptstyle #4}
!{(-3,-3)}="e"
!{(-3,-4)}*{\scriptscriptstyle #2}
!{(1,-3)}="f"
!{(1,-4)}*{\scriptscriptstyle #3}
"a"-"b"
"b"-"c"
"b"-"d"
"c"-"e"
"c"-"f"
}
}}

\newcommand{\ZX}[6][{}]{\ensuremath{
 \xygraph{
!{<0pt,0pt>;<4pt,0pt>:<0pt,-4pt>::}
!{(1,4)}*{\scriptscriptstyle #1}
!{(1,3)}="a"
!{(1,1)}*{\scriptscriptstyle{\circ}}="b"
!{(2,1)}*{\scriptscriptstyle #6}
!{(-1,-1)}*{\scriptscriptstyle{\bullet}}="c"
!{(0,-1)}*{\scriptscriptstyle #5}
!{(3,-1)}="d"
!{(3,-2)}*{\scriptscriptstyle #4}
!{(-3,-3)}="e"
!{(-3,-4)}*{\scriptscriptstyle #2}
!{(1,-3)}="f"
!{(1,-4)}*{\scriptscriptstyle #3}
"a"-"b"
"b"-"c"
"b"-"d"
"c"-"e"
"c"-"f"
}
}}

\newcommand{\Xsmall}[1]{\ensuremath{ 
 \xygraph{
!{<0pt,0pt>;<2.6pt,0pt>:<0pt,-2.6pt>::}
!{(0,2)}="a"
!{(0,0)}*{\scriptscriptstyle{\circ}}="b"
!{(1.4,0)}*{\scriptscriptstyle{#1}}
!{(-2,-2)}="c"
!{(2,-2)}="d"
"a"-"b"
"b"-"c"
"b"-"d"
}
}}

\newcommand{\Ysmall}[1]{\ensuremath{ 
 \xygraph{
!{<0pt,0pt>;<2.6pt,0pt>:<0pt,-2.6pt>::}
!{(0,2)}="a"
!{(0,0)}="b"
!{(1.4,0)}*{\scriptscriptstyle{#1}}
!{(-2,-2)}="c"
!{(2,-2)}="d"
"a"-"b"
"b"-"c"
"b"-"d"
}
}}

\newcommand{\Zsmall}[1]{\ensuremath{ 
 \xygraph{
!{<0pt,0pt>;<2.6pt,0pt>:<0pt,-2.6pt>::}
!{(0,2)}="a"
!{(0,0)}*{\scriptscriptstyle{\bullet}}="b"
!{(1.4,0)}*{\scriptscriptstyle{#1}}
!{(-2,-2)}="c"
!{(2,-2)}="d"
"a"-"b"
"b"-"c"
"b"-"d"
}
}}

\begin{document}

\title[Operads of compatible structures and weighted partitions]{Operads of compatible structures and weighted partitions}

\author{Henrik Strohmayer}
\address{Department of Mathematics \\ Stockholm University \\ 106 91 Stockholm \\ Sweden}
\email{henriks@math.su.se}

\begin{abstract}
In this paper we describe operads encoding two different kinds of compatibility of algebraic structures. We show that there exist decompositions of these in terms of black and white products and we prove that they are Koszul for a large class of algebraic structures by using the poset method of B.~Vallette. In particular we show that this is true for the operads of compatible Lie, associative and pre-Lie algebras.
 \end{abstract}

\maketitle

\section*{Introduction}
\label{compatibility}

Let $[\shs\circ\shs]$ and $[\shs\bullet\shs]$ be Lie brackets on a common vector space over a field $\K$. One can then define a new bracket $[\shs ,\shs]$ by $[a,b]:=\alpha[a\circ b]+\beta[a\bullet b]$, for some $\alpha,\beta\in\K$. Any such bracket is clearly skew symmetric and bilinear, so the only condition necessary in order for $[\shs,\shs]$ to be a Lie bracket is that the Jacobi identity $[[a,b],c]+[[b,c],a]+[[c,a],a]=0$ should hold. Direct calculation shows that this condition is equivalent to
\[[[a\circ b] \bullet c]+[[b\circ c]\bullet a]+[[c\circ a]\bullet b]+[[a\bullet b]\circ c]+[[b\bullet c]\circ a]+[[c\bullet a]\circ b]=0 .\] 

This notion of compatibility of Lie algebras was considered already in \cite{Magri1978}, Equation (3.1), in the study of integrable Hamiltonian equations. There the condition was considered for two symplectic operators and F.~Magri called it the coupling condition.

It is the aim of this paper to study this kind of compatibility for any algebraic structure given by a binary quadratic operad.

\begin{defnalpha}
 Let $\Ocal$ be a binary quadratic operad and $U$ a vector space over $\K$. Let $A=(U,\mu_1\dotsc ,\mu_k)$ and $B=(U,\nu_1\dotsc ,\nu_k)$ be $\Ocal$-algebra structures on $U$. Define new operations by $\eta_i:=\alpha\mu_i+\beta\nu_i$ for some $\alpha,\beta\in\K$. We say that $A$ and $B$ are \emph{linearly compatible} if $C=(U,\eta_1,\dotsc ,\eta_k)$ is an $\Ocal$-algebra for any choice of $\alpha$ and $\beta$. Note that this is equivalent to requiring $C$ to be an $\Ocal$-algebra for $\alpha=\beta=1$.
\end{defnalpha}

In \cite{Dotsenko2007}, A.~Khoroshkin and V.~Dotsenko described the operad $\Lietwo$ encoding two compatible Lie algebras. They also considered the Koszul dual operad $\Comtwo$ encoding two compatible commutative algebras. The compatibility condition in this case is quite different from the linear compatibility of Lie algebras. The commutative associative products $\circ$ and $\bullet$ are compatible in the sense that, firstly it should not matter in which order the products appear, i.e.~that
\[
(a\circ b) \bullet c = (a\bullet b) \circ c , 
\]
and secondly that the associativity relation should be fulfilled for any order of applying the two products, i.e.
\[
(a\circ b) \bullet c = a\bullet (b \circ c),\bhs  (a\bullet b)\circ c = a\circ(b \bullet c).
\]

Structures compatible in this way we call \emph{totally compatible} since $\circ$ and $\bullet$ are totally interchangeable up to the number of each of them. 
We will give a formal definition of the notion of total compatibility by defining the related operad in Section \ref{compatibilityencodedbyoperads}.

In \cite{Vallette2007}~B.~Vallette introduced a new method for showing the Koszulness of algebraic operads which can be obtained as the linearization of a set operad. By associating a certain poset to a set operad $\Pcal$, and then studying its Cohen-Macaulay properties, one gets a concrete recipe for checking whether the algebraic operad associated to $\Pcal$, and thus also its Koszul dual operad, is Koszul or not. Studying the posets of unordered and ordered pointed and multipointed partitions in \cite{Chapoton2006}, B.~Vallette and F.~Chapoton were able to prove the Koszulness of several important operads such as $\Perm$, $\Prelie$, $\Comtrias$, $\Postlie$, $\Dias$, $\Dend$, $\Trias$ and $\Tridend$ over a field of any characteristic and over $\Z$. 

To show the Koszulness of $\Lietwo$ and $\Comtwo$, as well as several other linearly and totally compatible structures, we will use the poset method of B.~Vallette. In order to handle the poset associated to an operad of two totally compatible structures we will show that it decomposes into the fiber product of two posets. The first one being the poset associated to the original structure and the other one being what we will call the poset of weighted partitions. In contrast to the posets studied by B.~Vallette and F.~Chapoton, these products of posets are not totally semimodular, therefore we need to refine the arguments of \cite{Chapoton2006} in order to show that they are Cohen-Macaulay.

The paper is organized as follows. In Section \ref{section1} we give some basic definitions about operads and recall the definitions of  $\Lietwo$ and $\Comtwo$. Then we give a description of the operads encoding compatible structures. We also show that there exist decompositions of the operads of compatible structures using black, white and Hadamard products. In Section \ref{section2} we define set operads and the posets associated to them. Thereafter we give the definition of the fiber product of posets and make some observations about how it relates to the Hadamard product of operads. In Section \ref{section3} we show how to describe the poset associated to $\Comtwo$ as the poset of weighted partitions, and then we proceed to prove the Koszulness of a class of operads of compatible structures.

All vector spaces and tensor products are considered over an arbitrary field $\K$. Given a finite set $S$ we denote its cardinality by $|S|$. By $\N$ we mean the set $\{1,2,\dotsc \}$. For $n\in\N$, we denote by $[n]$ the set $\{1,\dotsc,n\}$.  Let $\s_n$ denote the symmetric group of permutations of $[n]$. By $\bbone_n$ we denote the trivial representation of $\s_n$ and by $\sgn_n$ the sign representation.

\section{Compatibility of structures encoded by operads}
\label{section1}

\subsection{Algebraic operads}

To fix the notation we start by giving some definitions concerning operads. For an introduction to operads see e.g.~\cite{Loday1996,Markl2002}.

\begin{defn}
An $\s_n$-\emph{module} is a vector space $V$ with a right action of $\s_n$. A collection $(V(n))_{n\in\N}$ of $\s_n$-modules is called an $\s$-\emph{module}.
\end{defn}

Define a monoidal product in the category of $\s$-modules by:
\[
{V\circ W}(n) = \bigoplus_{1\leq k \leq n} \left(\bigoplus_{i_1+\dotsb+i_k=n} V_k \otimes (W_{i_1} \otimes \dotsb \otimes W_{i_k} ) 
\otimes_{\s_{i_1}\times\dotsb\times\s_{i_k}}\s_n \right)_{\s_k},
\]
where we consider the coinvariants with respect to the action of $\s_k$ given by 
$(v\otimes(w_{i_1}\otimes\dotsb\otimes w_{i_k})\otimes\sigma) \tau=(v \tau \otimes (w_{i_{\tau(1)}}\otimes\dotsb \otimes w_{i_{\tau(k)}} )\otimes\bar{\tau}^{-1} \sigma)$ and $\bar{\tau}$ is the induced block permutation. A unit $I$ with respect to this product is given by the $\s$-module defined by
\[
I_n:= \left\{ \begin{array}{cl}
               \bbone_1       & \mbox{if $n=1$}   \\
               0              & \mbox{if $n\neq 1$}  
              \end{array}
       \right. .
\]

\begin{defn}
An \emph{algebraic operad} is a monoid $(\Ocal,\mu\colon\Ocal\circ\Ocal\to\Ocal,\Ecal\colon I \to\Ocal)$ in the monoidal category $(\smodules,\circ, I)$. For an element $(e\otimes(e_1 \otimes \dotsb \otimes e_k)\otimes \sigma)\in\Ocal\circ\Ocal$ we will suppress the sigma and denote $\mu(e\otimes(e_1\otimes\dotsb\otimes e_k))$ by $\mu(e;e_1,\dotsc, e_k)$.
\end{defn}

The fundamental example is the endomorphism operad.

\begin{ex}
 For a vector space $U$ we let $\Ecal_U(n):=\Hom_{\K}(U^{\otimes n} ,U)$ and define $\mu(f;f_1,\dotsc , f_k)$ to be the usual composition of multivariable functions. $\Ecal_U:=(\Ecal_U(n))_{n\in\N}$ is then an operad.
\end{ex}

\begin{defn}
 A \emph{representation} of an operad $\Ocal$ in a vector space $U$ is a morphism of operads $\rho:\Ocal\to\Ecal_U$. We call a vector space equipped with this extra structure an $\Ocal$-\emph{algebra}.
\end{defn}

The free operad on an $\smodule$ $V$, $\Free(V)$ can be described as follows. The part $\Free(V)(n)$ is freely generated by all rooted directed trees with $n$ leaves whose internal vertices are decorated by elements of $V$ and whose leaves are labeled by $[n]$. We consider the direction to be towards the root, thus an internal vertex $v$ always has one outgoing edge and one or more incoming edges. We decorate an internal vertex with an element of $V(m)$ if it has $m$ incoming edges. The composition product $\mu(T;T_1,\dotsc,T_k)$ is given by grafting the roots of $(T_1,\dotsc ,T_k)$ to the $k$ leafs of $T$. We let $\Free_{(i)}(V)$ denote all labeled decorated trees with $i$ internal vertices.

\begin{defn}
A \emph{quadratic operad} $\QO{V}{R}$ is the free operad on an $\smodule$ $V$ modulo relations $R\subset \Free_{(2)}(V)$. We call a quadratic operad \emph{binary} if $V(n)=0$ for all $n\neq 2$.
\end{defn}

\begin{rem}
 Let $\Ocal=\QO{V}{R}$ be a binary quadratic operad with a $\K$-basis $e_1,\dotsc , e_s$ of $V$. A representation $\rho$ of $\Ocal$ in a vector space $U$ can be thought of as the data $(U,\{\rho(e_1),\dotsc,\rho(e_s)\})$, where the $\rho(e_i)$ are binary operations on $U$ subject to axioms encoded by the relations $R$ and the $\smodule$ strucure of $V$.
\end{rem}

For an $\smodule$ $V$ concentrated in $V(2)$ we have that $\Free_{(2)}(V)=\Free_{(2)}(V)(3)$. Given such a module $V$ with $\K$-basis  $\{\Ysmall{1},\dotsc, \Ysmall{s}\}$ we denote a labeled tree in $\Free_{(2)}(V)$ decorated with $\Ysmall{i}$ above $\Ysmall{j}$ by 
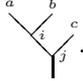
\[
 \YY{a}{b}{c}{i}{j}.
\]
The space $\Free_{(2)}(V)$ is then spanned by the trees
\[
\left\{\YY{1}{2}{3}{i}{j}, \YY{2}{3}{1}{i}{j}, \YY{3}{1}{2}{i}{j}\right\}_{1\leq i,j\leq s}.
\]
and is $3s^2$-dimensional. Thus for a binary quadratic operad $\Ocal=\QO{V}{R}$ operad we have that $R$ consists of $t\leq 3s^2$ linearly independent relations 
\begin{equation}
\label{quadraticrelations}
 R=\left\{\sum_{1\leq i,j \leq s} \gamma^{k,1}_{i,j}\YY{1}{2}{3}{i}{j}+\gamma^{k,2}_{i,j}\YY{2}{3}{1}{i}{j}+\gamma^{k,3}_{i,j}\YY{3}{1}{2}{i}{j} \right\}_{1\leq k \leq t}.
\end{equation}

If for an $\smodule$ $V$ it is true that $V(n)$ is finite dimensional for all $n$, then we can define its \emph{Czech dual} $V^{\cz}:=V(n)^*\otimes\sgn_n$. There is a natural pairing with respect to this duality,  $\langle \_ ,\_ \rangle \colon \Free_{(2)}(V^{\cz})\otimes \Free_{(2)}(V) \to \K$, given by
\[
 \left( \YY{a}{b}{c}{\hspace{7pt}i^{\cz}}{\hspace{7pt}j^{\cz}}, \YY{d}{e}{f}{k}{l}\right)=\delta_{(a,b,c),(d,e,f)}\delta_{i,k}\delta_{j,l}.
\]
We denote by $R^\perp$ the relations orthogonal to $R$ with respect to this pairing.

\begin{defn}
Let  $\Ocal=\QO{V}{R}$  be a quadratic operad. If $V(n)$ is finite dimensional for all $n$, then its \emph{Koszul dual operad} is given by $\Ocal^{!}:=\QO{V^{\cz}}{R^\perp}$.
\end{defn}

Given relations $R$ as in \eqref{quadraticrelations}, there are $3s^2-t$ linearly independent orthogonal relations 
\begin{equation}
\label{orthogonalrelations}
 R^\perp= \left\{\sum_{1\leq i,j \leq s} \eta^{k,1}_{i,j}\XX{1}{2}{3}{\hspace{7pt}i^{\cz}}{\hspace{7pt}j^{\cz}}+\eta^{k,2}_{i,j}\XX{2}{3}{1}{\hspace{7pt}i^{\cz}}{\hspace{7pt}j^{\cz}}+\eta^{k,3}_{i,j}\XX{3}{1}{2}{\hspace{7pt}i^{\cz}}{\hspace{7pt}j^{\cz}}\right\}_{1\leq k \leq 3s^2-t}.
\end{equation}

\begin{rem}
For a binary quadratic operad  $\Ocal=\QO{V}{R}$ such that $V(n)$ is finite dimensional for all $n\in\N$ we have that $(\Ocal^!)^!=\Ocal$.
\end{rem}

\subsection{The operads $\Comtwo$ and $\Lietwo$}
\label{lietwocomtwo}

\begin{defn}
The operad $\Lietwo$, encoding two linearly compatible Lie algebras, is the quadratic operad $\QO{V}{R}$ where the $\smodule$ $V$ is given by
\[
V(n):= \left\{ \begin{array}{ll}
               \sgn_2\oplus \sgn_2 & \mbox{if $n=2$}   \\
               0                   & \mbox{if $n\neq 2$} 
              \end{array}
       \right. .
\]
We represent a natural $\K$-basis of $V(2)$ as two binary corollas
\[
\sgn_2\oplus \sgn_2=\K\W{1}{2}\oplus\K\B{1}{2}.
\]
Then the relations $R$ are as follows
\begin{align*}
&\WW{1}{2}{3}+\WW{2}{3}{1}+\WW{3}{1}{2}, \bhs \BB{1}{2}{3}+\BB{2}{3}{1}+\BB{3}{1}{2}, \\
&\WB{1}{2}{3}+\WB{2}{3}{1}+\WB{3}{1}{2}+\BW{1}{2}{3}+\BW{2}{3}{1}+\BW{3}{1}{2}.
\end{align*}
 
The Koszul dual operad $(\Lietwo)^{!}$ is generated by the $\smodule$ $V^{\cz}$ which has as only non-zero part
\[
V^{\cz}(2)=\bbone_2\oplus\bbone_2=\K\Wczech{1}{2}\oplus\K\Bczech{1}{2}.
\]
From now on we will omit $\vphantom{1}^{\cz}$ from the notation. The relations $R^{\perp}$ are then given by
\begin{align*}
&\WW{1}{2}{3}-\WW{2}{3}{1}, \bhs \WW{1}{2}{3}-\WW{3}{1}{2}, \bhs \BB{1}{2}{3}-\BB{2}{3}{1}, \bhs \BB{1}{2}{3}-\BB{3}{1}{2},\\
&\WB{1}{2}{3}-\BW{1}{2}{3},\bhs \WB{2}{3}{1}-\BW{2}{3}{1}, \bhs \WB{3}{1}{2}-\BW{3}{1}{2},\\
&\WB{1}{2}{3}-\WB{2}{3}{1},\bhs \WB{1}{2}{3}-\WB{3}{1}{2}, \bhs \BW{1}{2}{3}-\BW{2}{3}{1},\bhs \BW{1}{2}{3}-\BW{3}{1}{2}. 
\end{align*}
Note that the last two relations are a consequence of the previous five.

We denote $(\Lietwo)^{!}$ by $\Comtwo$. 
\end{defn}

\begin{rem}
The operads $\Lietwo$ and $\Comtwo$ were denoted by $\Lie_2$ and $\Com_2$, respectively, in \cite{Dotsenko2007}. We denote them with upper index so that the notation does not interfere with that of set operads, see Section \ref{setoperads}, and we have the indices to the right and left respectively to emphasize that these are two different kinds of compatibility.
\end{rem}

\begin{prop}
\label{comtwobasisprop}
We have that $\Comtwo(n)=\bbone_n\oplus\dotsb\oplus\bbone_n$, where the sum consists of $n$ terms.  In terms of labeled trees decorated with $V^{\cz}$, a $\K$-basis for $\Comtwo(n)$ is given by
\[
\left\{ 
\WBcorolla 
\right\}_{0\leq i \leq n-1} .
\] 

Denote by $D_{i}^n$ the basis element in $\Comtwo(n)$ corresponding to $i$ white products. The composition product in $\Comtwo$ is then given by
\[
\mu(D^n_{i};D_{i_1}^{m_1},\dotsc, D_{i_n}^{m_n}) = D^{m_1+\dotsb+m_n}_{i+i_1+\dotsb+i_n}
\]
\end{prop}
\begin{proof}
This follows from the relations of $\Comtwo$ being homogenous in the number of black and white products. Thus any element of $\Comtwo(n)$ is determined by the number of white products, which can be at most $n-1$.
\end{proof}

\subsection{Definitions of compatible structures}
\label{compatibilityencodedbyoperads}

Let $\Ocal=\QO{V}{R}$ be a binary quadratic operad. Consider two operads $\Ocal_\circ=\QO{V_\circ}{R_\circ}$ and $\Ocal_\bullet=\QO{V_\bullet}{R_\bullet}$ both isomorphic to $\Ocal$. We choose $\K$-bases $\Xsmall{1},\dotsc,\Xsmall{s}$ of $V_\circ$ and $\Zsmall{1},\dotsc,\Zsmall{s}$ of $V_\bullet$ so that there exists an isomorphism $\phi \colon\Ocal_\circ\to\Ocal_\bullet$ with $\phi(\Xsmall{i})=\Zsmall{i}$. The relations $R_\circ$ and $R_\bullet$ can then be given by the same $\gamma^{k,l}_{i,j}$, cf.~\eqref{quadraticrelations}. By embedding $\Ocal_\circ$ and $\Ocal_\bullet$ into $\QO{V_\circ\oplus V_\bullet}{R_\circ\cup R_\bullet}$ we obtain an operad whose representations are pairs of $\Ocal$-algebras which not necessarily are compatible in any way. In order to encode linear compatibility we define the following relations.
\begin{align*}
 R_{\circ\bullet}:=\left\{\sum_{1\leq i,j \leq s}\right. &\left.\gamma^{k,1}_{i,j}\XZ{1}{2}{3}{i}{j}+\gamma^{k,2}_{i,j}\XZ{2}{3}{1}{i}{j}+\gamma^{k,3}_{i,j}\XZ{3}{1}{2}{i}{j}+\right.\\
 &\left.\gamma^{k,1}_{i,j}\ZX{1}{2}{3}{i}{j}+\gamma^{k,2}_{i,j}\ZX{2}{3}{1}{i}{j}+\gamma^{k,3}_{i,j}\ZX{3}{1}{2}{i}{j}\right\}_{1\leq k \leq t}.
\end{align*}

\begin{defn}
Given $\Ocal_\circ$ and $\Ocal_\bullet$ and $R_{\circ\bullet}$ as above we define a new binary quadratic operad by $\lc{\Ocal}:=\QO{V_\circ\oplus V_\bullet}{R_\circ\cup R_\bullet \cup R_{\circ\bullet}}$.
\end{defn}

\begin{prop}
\label{lcprop}
 A representation of $\lc{\Ocal}$ is a pair of linearly compatible $\Ocal$-algebras.
\end{prop}
\begin{proof}
By direct calculation.
\end{proof}

We now turn our attention to the other kind of compatibility which should generalize the compatibility of  $\Comtwo$. Given a binary quadratic operad $\Ocal$ and isomorphic operads $\Ocal_\circ$ and $\Ocal_\bullet$ as above, we define
\[
 \mbox{$_{\circ\bullet}^{\phantom{\circ}1} R$} := \left\{\XZ{1}{2}{3}{i}{j}-\ZX{1}{2}{3}{i}{j},\XZ{2}{3}{1}{i}{j}-\ZX{2}{3}{1}{i}{j},\XZ{3}{1}{2}{i}{j}-\ZX{3}{1}{2}{i}{j}\right\}_{1\leq i,j \leq s}.
\]
These relations encode that the order in which we apply operations of $\Ocal_\circ$ and $\Ocal_\bullet$ is irrelevant. Next we define
\[
 \mbox{$_{\circ\bullet}^{\phantom{\circ}2} R$} :=
\left\{\sum_{1\leq i,j \leq s}  \gamma^{k,1}_{i,j}\XZ{1}{2}{3}{i}{j}+\gamma^{k,2}_{i,j}\XZ{2}{3}{1}{i}{j}+\gamma^{k,3}_{i,j}\XZ{3}{1}{2}{i}{j}\right\}_{1\leq k \leq t},
\]
which encode that the relations of the original operad is satisfied for combinations of the operations of $\Ocal_\circ$ and $\Ocal_\bullet$ if we first apply an operation of $\Ocal_\circ$ and then an operation of $\Ocal_\bullet$. Note that a consequence of  $\mbox{$_{\circ\bullet}^{\phantom{\circ}1} R$}$ and $\mbox{$_{\circ\bullet}^{\phantom{\circ}2} R$}$ is 
\[
 \mbox{$_{\bullet\circ}^{\phantom{\circ}2} R$} :=
\left\{\sum_{1\leq i,j \leq s}  \gamma^{k,1}_{i,j}\ZX{1}{2}{3}{i}{j}+\gamma^{k,2}_{i,j}\ZX{2}{3}{1}{i}{j}+\gamma^{k,3}_{i,j}\ZX{3}{1}{2}{i}{j}\right\}_{1\leq k \leq t}.
\]
We define  $ \mbox{$_{\circ\bullet} R$}:= \mbox{$_{\circ\bullet}^{\phantom{\circ}1} R$} \cup \mbox{$_{\circ\bullet}^{\phantom{\circ}2} R$}$.
\begin{defn}
Given $\Ocal_\circ$, $\Ocal_\bullet$ and $\mbox{$_{\circ\bullet} R$}$ as above we define a new binary quadratic operad by $\tc{\Ocal}:=\QO{V_\circ\oplus V_\bullet}{R_\circ\cup R_\bullet \cup \mbox{$_{\circ\bullet} R$}}$. A representation of $\tc{\Ocal}$ is a pair of $\Ocal$-algebras with the compatibility given by $\mbox{$_{\circ\bullet} R$}$.  We call structures compatible in this way \emph{totally compatible}.
\end{defn}

We note that $\Comtwo$ is an operad of this form.

\begin{prop}
\label{lckoszuldualtotc}
Let $\Ocal=\QO{V}{R}$ be a binary quadratic operad such that $V(n)$ is finite dimensional for all $n\in\N$. We have $(\lc{\Ocal})^!=\tc{(\Ocal^!)}$ and $(\tc{\Ocal})^!=\lc{(\Ocal^!)}$.
 \end{prop}
\begin{proof}
By direct calculation.
\end{proof}

\subsection{Black product, white product and Hadamard product}

In \cite{Ginzburg1994,Ginzburg1995} V.~Ginzburg and M.~Kapranov generalized the notions of black and white products for algebras to binary quadratic operads. B.~Vallette generalized the notion further to arbitrary operads given by generators and relations and to properads in \cite{Vallette2006}. He also established some results about black and white products for operads. For more details we refer the reader to these papers.

The definition of the black product for binary quadratic operads is given in terms of a certain map $\Psi$. Note that for binary quadratic operads $\Free(V)(3)$ is equal to $\Free_{(2)}(V)$ and that $\Free_{(2)}(V)$ is spanned by three types of decorated trees, corresponding to the possible labelings of the leaves, see Section \ref{compatibilityencodedbyoperads}. Given two binary quadratic operads $\Ocal=\QO{V}{R}$ and $\Qcal=\QO{W}{S}$ the map
\[
 \Psi\colon\Free(V)(3)\otimes \Free(W)(3) \otimes \sgn_3 \to \Free(V\otimes W \otimes \sgn_2)
\]
is defined by 
\[
 \YY{a}{b}{c}{i}{j}\otimes\YY{d}{e}{f}{k}{l}\mapsto \delta_{(a,b,c),(d,e,f)} \YY{a}{b}{c}{\hspace{8pt} i\otimes k}{\hspace{8pt}j\otimes l}\shs,
\]
where by abuse of notation $i\otimes k$ denotes the tensor product of the elements decorating the trees.

\begin{defn}
\label{blackproduct}
 Let $\Ocal=\QO{V}{R}$ and $\Qcal=\QO{W}{S}$ be binary quadratic operads whose $\smodules$ of generators $V$ and $W$ are finite dimensional. We define their \emph{black product} by
\[
 \Ocal\bullet\Qcal:=\QO{V\otimes W \otimes \sgn_2}{\Psi(R\otimes S)}.
\]
\end{defn}

The white product is defined through another map
\[
 \Phi \colon  \Free(V\otimes W)(3) \to \Free(V)(3)\otimes \Free(W)(3)
\]
which is given by
\[
 \YY{a}{b}{c}{\hspace{9pt} i\otimes k}{\hspace{9pt}j\otimes l} \mapsto \YY{a}{b}{c}{\hspace{2pt}i}{\hspace{2pt}j}\otimes\YY{a}{b}{c}{\hspace{2pt} k}{\hspace{2pt}l}.
\]
\begin{defn}
  Let $\Ocal=\QO{V}{R}$ and $\Qcal=\QO{W}{S}$ be binary quadratic operads. We define their \emph{white product} by
\[
 \Ocal\circ\Qcal:=\QO{V\otimes W}{\Phi^{-1}(R\otimes\Free(W)(3)+\Free(V)(3)\otimes S)}.
\]
\end{defn}

We have the following relation between black and white products which was stated in \cite{Ginzburg1994,Ginzburg1995} and explicitly proven in \cite{Vallette2006}.

\begin{prop}[Theorem 2.2.6 in \cite{Ginzburg1994}]
\label{blackwhitekoszuldual}
 Let $\Ocal$ and $\Qcal$ be binary quadratic operads generated by finite dimensional $\smodules$, then $(\Ocal\circ\Qcal)^!=\Ocal^!\bullet\Qcal^{!}$.
\end{prop}

We now reach the highlight of this section with the following theorem.

\begin{thm}
\label{blackproductlinearcompatibility}
 Let $\Ocal$ be a binary quadratic operad. We have $\lc{\Ocal}=\Ocal\bullet\Lietwo$.
\end{thm}
\begin{proof}
 Let $\Ocal=\QO{V}{R}$ be generated by the $\smodule$ $V=\K \Ysmall{1}\oplus\dotsb\oplus \K \Ysmall{s}$ and with relations
\[
R=\left\{\sum_{1\leq i,j \leq s} \gamma^{k,1}_{i,j}\YY{1}{2}{3}{i}{j}+\gamma^{k,2}_{i,j}\YY{2}{3}{1}{i}{j}+\gamma^{k,3}_{i,j}\YY{3}{1}{2}{i}{j} \right\}_{1\leq k \leq t}.
\]
 Further denote the generators of $\Lietwo$ as in Section \ref{lietwocomtwo} by $\Wsmall$ and $\Bsmall$ and the relations by $S=S_\circ\cup S_\bullet\cup S_{\circ\bullet}$. 

By Definition \ref{blackproduct} we see that $\Ocal\bullet\Lietwo$ is generated by $(\K \Ysmall{1}\oplus\dotsb\oplus \K \Ysmall{s}) \otimes (\K\Wsmall\oplus\K\Bsmall) \otimes \sgn_2$. We denote the generator of $\K\Ysmall{i}\otimes\K\Wsmall\otimes\sgn_2$ by $\Xsmall{i}$ and since $\K\Wsmall\otimes\sgn_2=\sgn_2\otimes\sgn_2\iso\bbone_2$ we have that $V_\circ= \K \Xsmall{1}\oplus\dotsb\oplus \K \Xsmall{s}$ is isomorphic to $V$ as an $\smodule$. Of course the same is true for  $V_\bullet=\K \Zsmall{1}\oplus\dotsb\oplus \K \Zsmall{s}$, with the obvious meaning of $\Zsmall{i}$.

Next we see that 
\[
 R_\circ=\Psi(R\otimes S_\circ)=\left\{\sum_{1\leq i,j \leq s} \gamma^{k,1}_{i,j}\XX{1}{2}{3}{i}{j}+\gamma^{k,2}_{i,j}\XX{2}{3}{1}{i}{j}+\gamma^{k,3}_{i,j}\XX{3}{1}{2}{i}{j} \right\}_{1\leq k \leq t}.
\]
and similarly for $R_\bullet=\Psi(R\otimes S_\bullet)$. Finally for $R_{\circ\bullet}=\Psi(R\otimes S_{\circ\bullet})$ we have
\begin{align*}
 R_{\circ\bullet}=\left\{\sum_{1\leq i,j \leq s}\right.   
&\left.\gamma^{k,1}_{i,j}\XZ{1}{2}{3}{i}{j}+\gamma^{k,2}_{i,j}\XZ{2}{3}{1}{i}{j}+\gamma^{k,3}_{i,j}\XZ{3}{1}{2}{i}{j}+\right.\\
&\left.\gamma^{k,1}_{i,j}\ZX{1}{2}{3}{i}{j}+\gamma^{k,2}_{i,j}\ZX{2}{3}{1}{i}{j}+\gamma^{k,3}_{i,j}\ZX{3}{1}{2}{i}{j}\right\}_{1\leq k \leq t}.
\end{align*}
where ${1\leq i \leq n}$. Thus we see that $\Ocal\bullet\Lietwo=\QO{V_\circ\oplus V_\bullet}{R_\circ\cup R_\bullet\cup R_{\circ\bullet}}=\lc{\Ocal}$.
\end{proof}

\begin{cor}
\label{whiteproducttotalcompatibility}
  Let $\Ocal$ be a binary quadratic operad. We have $\tc{\Ocal}=\Ocal\circ\Comtwo$.
\end{cor}
\begin{proof}
 By Proposition \ref{lckoszuldualtotc}  we have that $(\lc{(\Ocal^!)})^!=\tc{\Ocal}$ and by Theorem \ref{blackproductlinearcompatibility} that  $\lc{(\Ocal^!)}=\Ocal^!\bullet\Lietwo$.  By Proposition \ref{blackwhitekoszuldual} we know that $(\Ocal^!\bullet\Lietwo)^!=\Ocal\circ\Comtwo$. Putting this together we conclude that 
\[\tc{\Ocal}=(\lc{(\Ocal^!)})^!=(\Ocal^!\bullet\Lietwo)^!=\Ocal\circ\Comtwo.
\]

\end{proof}

In practice the white product can be difficult to compute explicitly. In \cite{Vallette2006} a useful result was proven relating the white product and Hadamard product for certain operads. 
\begin{defn}

The \emph{Hadamard product} $\Ocal\otimes_H \Qcal$ of two operads $\Ocal$ and $\Qcal$ is defined as $(\Ocal\otimes_H \Qcal)(n):=\Ocal(n)\otimes \Qcal(n)$. The composition $\mu$ is given by
\[
\mu(e\otimes q;e_1\otimes q_1,\dotsc,e_k \otimes q_k):= \mu(e;e_1,,\dotsc,e_k) \otimes \mu(q;q_1,\dotsc,q_k).
\]
\end{defn}

For a quadratic operad $\Ocal=\QO{V}{R}$ let $\pi_\Ocal\colon\Free(V)\to\Ocal$ be the natural projection. Denote by $T$ a labeled binary tree with $n-1$ internal vertices. We order the internal vertices linearly in an arbitrary way and let $\Lcal^V_T$ denote the induced decoration morphism $\Lcal^V_T\colon V^{\otimes (n-1)}\to\Free(V)$ which decorates the internal vertices of $T$ with elements of $V$.

\begin{prop}[Proposition 15 in \cite{Vallette2007}]
\label{surjectivewhitehadamard}
 Let $\Ocal$ be a binary quadratic operad such that for every $n\geq 3$ and every labeled binary tree $T$ with $n-1$ vertices the composite map $\pi_\Ocal\circ\Lcal_T^V\colon V^{\otimes (n-1)}\to\Ocal(n)$ is surjective. For every binary quadratic operad, $\Qcal$, the white product $\Ocal\circ\Qcal$ is equal to the Hadamard product $\Ocal \otimes_H \Qcal$.
\end{prop}

Since we will use the condition in Proposition \ref{surjectivewhitehadamard} later we extract it into a definition.

\begin{defn}
\label{weakassdef}
Let $\Ocal=\QO{V}{R}$ be a binary quadratic operad with a $\K$-basis  $\{\Ysmall{1},\dotsc , \Ysmall{s}\}$ of $V$. Denote the element $\Ysmall{i}(12)$ by $\Ysmall{\hspace{7pt}i^{\op}}$. We call $\Ocal$ \emph{weakly associative} if
\[
 \forall \XX{a}{b}{c}{i}{j}  \shs \exists  \XX{b}{c}{a}{k}{\hspace{7pt}l^{\op}} \bhs \mbox{such that} \bhs \XX{a}{b}{c}{i}{j} = \XX{b}{c}{a}{k}{\hspace{7pt}l^{\op}}.
\]
\end{defn}

Note that an operation $\Ysmall{i}$ is associative in the usual sense precisely when the above condition is satisfied for $i=j=k=l$. 

\begin{prop}
\label{weakassprop}
 Let $\Ocal$ be binary quadratic operad, then $\Ocal$ is weakly associative iff it has the property of Proposition \ref{surjectivewhitehadamard}.
\end{prop}

\begin{proof}
 Assume that $\Ocal=\QO{V}{R}$ is weakly associative. Let $T$ be any labeled binary tree. By repeatedly using the identity
\[
 \XX{a}{b}{c}{i}{j} = \XX{b}{c}{a}{k}{\hspace{7pt}l^{\op}}
\]
 any decorated labeled binary tree $T'$ is equivalent to a decorated tree with the same shape and labeling as $T$. Hence the map $\pi_{\subtc{\Com}}\circ\Lcal_T^V$ is surjective.

 Now assume instead that $\pi_{\subtc{\Com}}\circ\Lcal_T^V$ is surjective for any labeled binary tree $T$. Let 
\[T=\XX{b}{c}{a}{{}}{{}}.\] 
Then since $\pi_{\subtc{\Com}}\circ\Lcal_T^V$ is surjective, any decorated tree
\[\XX{a}{b}{c}{i}{j}\]
is equivalent to a decorated tree of the same shape and labeling as $T$, which is exactly the condition in Definition \ref{weakassdef}.
\end{proof}

\begin{cor}
\label{comtwowhiteproducthadamardproduct}
 For every binary quadratic operad $\Ocal$ we have $\Ocal\circ\Comtwo=\Ocal\otimes_H\Comtwo$.
\end{cor}
\begin{proof}
Clearly $\Comtwo$ is weakly associative, thus by Proposition~\ref{weakassprop} it satisfies the condition of Proposition~\ref{surjectivewhitehadamard} whence we obtain the desired result.
\end{proof}
\section{Operadic partition posets of set operads}
\label{section2}

\subsection{Set operads}
\label{setoperads}

An $\s$-set is a collection of sets, $S=(S_n)_{n\in\N}$, equipped with a right action of the symmetric group $\s_{n}$ on $S_n$. Define a monoidal product in the category of $\smodules$ by:
\[
{S\circ T}_n = \bigsqcup_{1\leq k \leq n} \left(\bigsqcup_{i_1+\dots+i_k=n} S_k \times (T_{i_1} \times \dotsb \times T_{i_k} ) 
\times_{\s_{i_1}\times\dotsb\times\s_{i_k}}\s_n \right)_{\s_k},
\]
where we consider the coinvariants with respect to the action of $\s_k$ given by 
$(s,(t_{i_1},\dots,t_{i_k}),\sigma) \tau=(s \tau,(t_{i_{\tau(1)}},\dotsc,t_{i_{\tau(k)}} ),\bar{\tau}^{-1} \sigma)$ and $\bar{\tau}$ is the induced block permutation. A unit $I$ with respect to this product is given by the $\smodule$ defined by
\[
I_n:= \left\{ \begin{array}{cl}
               [1]       & \mbox{if $n=1$}   \\
               \emptyset & \mbox{if $n\neq 1$} . 
              \end{array}
       \right.
\]

\begin{defn}
A \emph{set operad} is a monoid $(\Pcal,\mu\colon\Pcal\circ\Pcal\to\Pcal,\varepsilon\colon I \to\Pcal)$ in the monoidal category $(\s\mbox{-sets},\circ, I)$. For an element $(p,(p_1,\dotsc, p_k),\sigma)\in\Pcal\circ\Pcal$ we will suppress the sigma and denote $\mu(p,(p_1,\dotsc, p_k))$ by $\mu(p;p_1,\dotsc, p_k)$.
\end{defn}

To any set operad $\Pcal$ one can associate an algebraic operad $\aao{\Pcal}$ by considering formal linear combinations of the elements, i.e.~ $\aao{\Pcal}(n)=\K[\Pcal_n]$. We call $\aao{\Pcal}$ the \emph{linearization} of $\Pcal$. Often we will use the same notation for a set operad as for its linearization. It should be clear from the context which of the two is referred to.

To an element $(p_1,\dotsc, p_k)\in\Pcal_{i_1}\times \dots \times \Pcal_{i_k}$ one can associate a map 
\[
\mu_{p_1,\dotsc, p_k} \colon \Pcal_{k}\to\Pcal_{i_1+\dots+i_k}
\]
defined as $\mu_{p_1,\dotsc, p_k}(p):=\mu(p ; p_1,\dotsc, p_k)$. The following definition was introduced in \cite{Vallette2007} since it is a crucial property for set operads in order to use the poset method.

\begin{defn}
A set operad $\Pcal$ is called a \emph{basic-set operad} if the map $\mu_{p_1,\dotsc, p_k}$ is injective 
 for all $(p_1,\dotsc,p_k)\in\Pcal_{i_1}\times \dots \times \Pcal_{i_k}$.
\end{defn}

\begin{prop}
\label{comtwobasicset}
The operad $\Comtwo$ is the linearization of a basic-set operad.
\end{prop}

\begin{proof}
The operad $\Comtwo$ is the linearization of $\Pcal$, where $\Pcal_n=\{D^n_i\}$ and the $D^n_i$ are as in Proposition \ref{comtwobasisprop}. That $\Pcal$ is basic-set is immediate from the formula for the composition product.
\end{proof}

\subsection{Operadic partition posets}
\label{operadicpartitionposets}

 For definitions of the various notions related to posets see \cite{Bjorner1983,Vallette2007}.

\begin{defn}
Let $\Pcal$ be a set operad. A $\Pcal$-\emph{partition} of $[n]$ is the data 
$\{(B_1,p_1),$ 
$\dotsc,(B_s,p_s)\}$, where $\{B_1,\dotsc,B_s\}$ is a partition of $[n]$ and $p_i\in\Pcal_{|B_i|}$. We let $\Pi_\Pcal(n)$ denote the set of all $\Pcal$-partitions of $[n]$ and let $\Pi_\Pcal$ denote the collection $\{\Pi_\Pcal(n)\}_{n\in\N}$. For an algebraic operad $\Ocal$ which is the linearization of a set operad $\Pcal$, i.e.~$\Ocal=\aao{\Pcal}$, we will sometimes write $\Pi_{\Ocal}$ for $\Pi_{\Pcal}$.
\end{defn}

\begin{rem}
\label{rem-partition}
One can think of this as enriching a partition with elements of an operad or, shifting the perspective, as labeling the input of the operation that an element $p_i\in\Pcal_{|B_i|}$ describes with the elements of the block $B_i$ instead of with $[|B_i|]$. E.g. one can identify
\[
\left( \{3,4,7\}, \WW{2}{3}{1} \right)\sim \WW{4}{7}{3}.
\]
The definition in \cite{Vallette2007}~uses ordered sequences of elements of the blocks instead of unordered blocks and then considers equivalence classes of pairs $(S_B,p)$, where $S_B$ is an ordered sequence of the elements of a block $B$ where each element appears exactly once and $p\in\Pcal_{|S_B|}$. E.g.
\[
\left( (3,4,7), \WW{2}{3}{1} \right) \sim \left( (4,7,3), \WW{1}{2}{3} \right) \sim \WW{4}{7}{3}.
\]
Our definition corresponds to choosing the representative of a class with the elements of the sequence in ascending order. In the following we will assume that, given a partition $\alpha=\{(A_1,p_1),\dotsc,(A_r,p_r)\}$, the elements of a block $A_i=\{a^i_1,\dotsc,a^i_{m_i}\}$ are indexed in ascending order, i.e.~$a^i_j<a^i_{j+1}$.
\end{rem}

Next we define a partial order on $\Pi_\Pcal(n)$ .
\begin{defn}
\label{operadicorder}
Let $\alpha=\{(A_1,p_1),\dotsc,(A_r,p_r)\}$ and $\beta=\{(B_1,q_1),\dotsc,(B_s,q_s)\}$ be two $\Pcal$-partitions of $[n]$.  We let $\alpha\leq\beta$ if
\begin{enumerate}
\item $\{A_1,\dotsc,A_r\}$ is a refinement of $\{B_1,\dotsc,B_s\}$, i.e.~each $B_j$ is the union of one or more $A_i$.

\item when $B_j=A_{i_1}\cup\dots\cup A_{i_t}$ then there exists a $p\in\Pcal_t$ such that $q_j=\mu(p;p_{i_1},\dotsc,p_{i_t})\sigma^{-1}$, where $\sigma\in\s_{|B_j|}$ is the obvious permutation associated to
\[
\left(
\begin{array}{l}
b^j_1 \dots  b^j_{|B_j|} \\ 
a^{i_1}_1 \dots a^{i_t}_{m_{i_t}} 
\end{array} 
\right).
\]
\end{enumerate}
We call $\Pi_\Pcal$ together with this partial order the \emph{operadic partition poset} of $\Pcal$. 
\end{defn}

\begin{rem}
We define the order in the opposite way to the one in \cite{Vallette2007} to make it correspond to the way it is defined in \cite{Chapoton2006}. Note that with this in mind our definition leads to the same ordering of the corresponding equivalence classes.
\end{rem}

\begin{ex}
Using the identification in Remark \ref{rem-partition} we see that in $\Pi_{\subtc{\Com}}(7)$
\[
 \left\{ \WW{1}{2}{6} ,\shs \Id{5} \shs , \BW{3}{4}{7}   \right\} \leq  \left\{ \WW{1}{2}{6} , \WWB{3}{4}{5}{7} \right\}
\]
since 
\[
\mu(\W{1}{2}; \shs \Id{5} \shs ,\BW{3}{4}{7} )=\WBWvar{5}{3}{4}{7}= \WWB{3}{4}{5}{7}\shs.
\]
\end{ex}

In \cite{Vallette2007}, Vallette  studied homological properties of the order complex associated to the partition poset of an operad. The following is the main result.

\begin{thm}[Theorem 9 of \cite{Vallette2007}]
\label{Bruno-main-thm}
Let $\Pcal$ be a basic-set quadratic operad. The operad $\aao{\Pcal}$ is Koszul iff each subposet $[\hat{0},\gamma]$ of each $\Pi_\Pcal(n)$ is Cohen-Macaulay, where $\gamma$ is a maximal element of $\Pi_\Pcal(n)$.
\end{thm}

\subsection{Fiber product of operadic partition posets}

In \cite{Bjorner2005} a product of posets was introduced under the name Segre product and a particular case studied. We prefer to call it fiber product because it corresponds to this categorical construction.

\begin{defn}
Let $P$,$Q$ and $S$ be posets. Given poset maps $f\colon P \to S$ and $g \colon Q \to S$ we define $P\times_{f,g}Q$, the \emph{fiber product} of $P$ and $Q$ over $f,g$, to be the subset of $P\times Q$ consisting of pairs $(p,q)$ such that $f(p)=g(q)$. The order on $P\times_{f,g} Q$ is induced by the order on $P\times Q$ which is given by $(p,q)\leq (p',q')$ if $p\leq p'$ and $q\leq q'$.
\end{defn}

Let $\Pi_n$ denote the poset of partitions of $[n]$ and let $\Pi$ denote the collection $\{\Pi_n\}_{n\in\N}$. Further, given operadic partition posets $\Pi_\Pcal$ and $\Pi_\Qcal$, let  $f\colon \Pi_\Pcal \to \Pi$ and $g \colon \Pi_\Qcal \to \Pi$ be the natural projections which sends an element $\alpha=\{(A_1,p_1),\dotsc,$ $(A_m,p_m) \}\in\Pi_\Pcal$ to the underlying partition $\{A_1,\dotsc,A_m \}$ and similarly for $g$. Then  $\Pi_\Pcal \times_{f,g} \Pi_\Qcal$ consists of pairs $(\alpha,\beta)$, where $\alpha=\{(A_1,p_1),\dotsc,(A_m,p_m) \}$ and $\beta=\{(A_1,q_1),\dotsc,$ $(A_m,q_m) \}$. This poset is isomorphic to the poset consisting of elements $\alpha=\{(A_1,p_1,q_1),\dotsc,(A_m,p_m,q_m)\}$, where $p_i\in\Pcal_{|A_i|}$ and $q_i\in\Qcal_{|A_i|}$, with the order given by $\alpha\leq\alpha'$ if
\begin{enumerate}
\item $\{A_1,\dotsc,A_r\}$ is a refinement of $\{A'_1,\dotsc,A'_s\}$.
\item when $A'_j=A_{i_1}\cup\dots\cup A_{i_t}$ then there exists a $p\in\Pcal_t$ and a $q\in\Qcal_t$ such that $p'_j=\mu(p;p_{i_1},\dotsc,p_{i_t})\sigma^{-1}$ and  $q'_j=\mu(q;q_{i_1},\dotsc,q_{i_t})\sigma^{-1}$, where $\sigma\in\s_{|A'_j|}$ is the permutation given in Definition~\ref{operadicorder}. 
\end{enumerate}
We will denote this fiber product by $\Pi_\Pcal \times_\Pi \Pi_\Qcal$. Note that $\Pi_{\Com}=\Pi$ whence $\Pi_\Pcal \times_\Pi \Pi_{\Com}= \Pi_\Pcal$, for any $\Pcal$.

\begin{defn}
The \emph{Hadamard product} $\Pcal\times_H \Qcal$ of two set operads $\Pcal$ and $\Qcal$ is defined as $(\Pcal\times_H \Qcal)_n=\Pcal_n\times \Qcal_n$, where $\times$ denotes the cartesian product. The composition $\mu$ is given by
\[
\mu((p,q);(p_1,q_1),\dotsc,(p_k,q_k)):=( \mu(p,;p_1,,\dotsc,p_k) ,\mu(q;q_1,\dotsc,q_k) ).
\]
\end{defn}

\begin{prop} 
\label{fiberproducthadamardproduct}
For any set operads $\Pcal$, $\Qcal$ the following equalities hold.
\begin{enumerate} 
 \item \label{fiberproducthadamardproductfirst}
 $\Pi_{\Pcal\times_H \Qcal}=\Pi_\Pcal \times_\Pi \Pi_\Qcal$
 \item \label{fiberproducthadamardproductsecond}
 $\aao{\Pcal\times_H \Qcal}=\aao{\Pcal}\otimes_H \aao{\Qcal}$
\end{enumerate} 
\end{prop}
\begin{proof}
 Immediate from the definitions involved.
\end{proof}

Next we describe the operadic partition poset associated to an operad encoding totally compatible structures.

\begin{cor}
\label{posetoftc}
Let $\Ocal$ be an algebraic operad which is the linearization of a set operad $\Pcal$. Then
\begin{enumerate}
\item $\tc{\Ocal}=\aao{\Pcal\times_H \Comtwo}$ and
\item $\Pi_{\subtc{\Ocal}}=\Pi_{\Pcal}\times_\Pi \Pi_{\subtc{\Com}}$.
\end{enumerate} 
\end{cor}

\begin{proof}
 Using  Corollary \ref{whiteproducttotalcompatibility} and Propositions \ref{comtwowhiteproducthadamardproduct} and \ref{fiberproducthadamardproduct} (\ref{fiberproducthadamardproductsecond}) we have that 
\[
 \tc{\Ocal}=\Ocal\circ\Comtwo=\Ocal\otimes_H\Comtwo=\aao{\Pcal\times_H \Comtwo}.
\]
Thus by Proposition \ref{fiberproducthadamardproduct} (\ref{fiberproducthadamardproductfirst}) we have $\Pi_{\subtc{\Ocal}}=\Pi_{\Pcal}\times_H \Pi_{\subtc{\Com}}$.
\end{proof}

We define $\tc{\Pcal}:=\Pcal\times_H \Comtwo$ and observe that $\aao{\tc{\Pcal}}=\tc{\aao{\Pcal}}$. 

\begin{prop}
\label{fiberproductbasicset}
Let $\Pcal$ and $\Qcal$ be set operads. If $\Pcal$ and $\Qcal$ are basic-set, then so is $\Pcal\times_H \Qcal$.
\end{prop}
\begin{proof}
 We want to show that the map $\mu_{(\nu_1,\eta_1),\dotsc,(\nu_k,\eta_k)}\colon\Pcal_k\times\Qcal_k\to\Pcal_{i_1+\dotsb+i_k}\times\Qcal_{i_1+\dotsb+i_k}$ given by $\mu_{(\nu_1,\eta_1),\dotsc,(\nu_k,\eta_k)}(\alpha,\beta)=(\mu(\alpha;\nu_1,\dotsc ,\nu_k),\mu(\beta;\eta_1,\dotsc,\eta_k))$ is injective for all $((\nu_1,\eta_1),\dotsc,(\nu_k,\eta_k))\in (\Pcal_{i_1}\times \Qcal_{i_1}) \times\dotsb\times (\Pcal_{i_k}\times \Qcal_{i_k})$. Now let $(\alpha,\beta),(\alpha',\beta')\in\Pcal_k\times\Qcal_k$ be such that $(\alpha,\beta)\neq(\alpha',\beta')$. Then $\alpha\neq\alpha'$ or $\beta\neq\beta'$ and thus, since $\Pcal$ and $\Qcal$ are basic set, either $\mu(\alpha;\nu_1,\dotsc ,\nu_k)\neq\mu(\alpha';\nu_1,\dotsc ,\nu_k)$ or $\mu(\beta;\eta_1,\dotsc,\eta_k))\neq\mu(\beta';\eta_1,\dotsc,\eta_k))$.
\end{proof}

\begin{cor}
\label{tcbasicset}
Let $\Pcal$ be a basic-set operad, then so is  $\tc{\Pcal}$.
\end{cor}
\begin{proof}
By Proposition \ref{comtwobasicset} we know that $\Comtwo$ is basic-set. Thus we can apply Proposition \ref{fiberproductbasicset} to $\tc{\Pcal}=\Pcal\times_H \Comtwo$.
\end{proof}
\section{Koszulness of a class of compatible structures}
\label{section3}

\subsection{$\Comtwo$ and weighted partitions}

Theorem \ref{Bruno-main-thm} was used in \cite{Vallette2006} and \cite{Chapoton2006} to show the Koszulness of several operads. There it was shown that for the associated posets all maximal intervals $[\hat{0},\gamma]$ were totally semimodular. Hence by Corollary 5.2 of \cite{Bjorner1983} they are CL-shellable and by Proposition 2.3 of the same paper shellable whence it follows that they are Cohen-Macaulay by Theorem 4.2 of \cite{Garsia1980}. The chain of implications is
\begin{equation}
\label{chainofimplications}
\text{totally semimodular} \implies \text{CL-shellable} \implies \text{shellable} \implies \text{Cohen-Macaulay}.
\end{equation}

\begin{defn}
A finite poset $P$ is called \emph{semimodular} if it is bounded, i.e.~ has a least and a greatest element, and for any distinct $\kappa,\lambda\in P$ covering a $\nu\in P$ there exists a $\omega\in P$ covering both $\kappa$ and $\lambda$. The poset $P$ is said to be \emph{totally semimodular} if it is bounded and all intervals $[\zeta,\xi]$ are semimodular.
\end{defn}

\begin{rem}
Contrary to the claims in \cite{Dotsenko2007}, the maximal chains of $\Pi_{\subtc{\Com}}$ are not necessarily totally semimodular. E.g.~consider the elements 
\[
(\W{1}{2},\Id{3},\Id{4}),(\Id{1},\Id{2},\W{3}{4})\in[(\Id{1},\Id{2},\Id{3},\Id{4}),(\WBB{1}{2}{3}{4})]\subset\Pi_{\subtc{\Com}}(4).
\]
They both cover $(\Id{1},\Id{2},\Id{3},\Id{4})$ but the only element covering both of them is $(\W{1}{2},\W{3}{4})$ which does not belong to the interval $[(\Id{1},\Id{2},\Id{3},\Id{4}),(\WBB{1}{2}{3}{4})]$.
\end{rem}

By the chain of implications (\ref{chainofimplications}) we see that to show Cohen-Macaulayness of $\Pi_{\subtc{\Com}}$ and thus Koszulness of $\Comtwo$, it is in fact sufficient to show that the maximal intervals of $\Pi_{\subtc{\Com}}$ are CL-shellable. A poset is CL-shellable if a certain kind of labeling of the maximal chains is possible, see \cite{Bjorner1983}. By Theorem 3.2 of \cite{Bjorner1983}, showing CL-shellability of a poset is equivalent to showing that it admits a recursive atom ordering. Recall that the atoms of a poset are the elements covering $\hat{0}$.

\begin{defn}
\label{recatomordering}
A graded poset $P$ admits a \emph{recursive atom ordering} if the length of the poset is 1 or if the length is greater than 1 and there is an ordering $\alpha_1,\dotsc,\alpha_m $ of the atoms of $P$ satisfying
\begin{enumerate}
\item \label{firstcriterion}
 For all $j\in[m]$, $[\alpha_j,\hat{1}]$ admits a recursive atom ordering in which the atoms of $[\alpha_j,\hat{1}]$ that come first in the ordering are those that cover some $\alpha_i$, where $i<j$.
\item \label{secondcriterion}
For all $i<j$, if $\alpha_i,\alpha_j<\lambda$ then there is a $k<j$, not necessarily distinct from $i$, and an element $\kappa\leq\lambda$ such that $\kappa$ covers both $\alpha_j$ and $\alpha_k$ 
\end{enumerate}
\end{defn}

We will soon see that $\Pi_{\subtc{\Com}}$ admits a recursive atom ordering, but first we make the structure of $\Pi_{\subtc{\Com}}$ explicit by the following partition poset.

\begin{defn}
Given a partition $\beta=\{B_1,\dotsc,B_s\}$ of $[n]$, we assign a weight $w_i\in\N$ to each block $B_i=\{b^i_1,\dotsc,b^i_{k_i} \}$, with  $0\leq w_i\leq k_{i} -1$. The weight of the block is denoted by $\weight(B_i):=w_i$. The weight of a partition $\beta$ is $\weight(\beta):=\weight(B_1)+\dotsb+\weight(B_s)$. A partition with this extra structure we call a \emph{weighted partition} and we denote the set of weighted partions of $[n]$ by $\wp{n}$. The collection $\{\wp{n}\}_{n\in\N}$ is denoted by $\wp{{}}$. 

Let $\nblocks(\beta)$ be the number of blocks of $\beta$. Then we can define a partial order on $\wp{n}$ by letting $\alpha\leq\beta$ if 
\begin{enumerate}
\item the partition of $\alpha$ is a refinement of the partition of $\beta$ and 
\item $\weight(\beta)-\weight(\alpha)\leq \nblocks(\alpha)-\nblocks(\beta)$.
\end{enumerate}

We call $\wp{{}}$ together with this partial order the \emph{poset of weighted partitions}.
\end{defn}

\begin{rem}
We see that the covering relation $\prec$ of the above partial order is given by
$\alpha\prec\beta$ if 
\begin{enumerate}
\item the partition of $\alpha$ is a refinement of that of $\beta$ obtained by splitting exactly one block of $\beta$ into two and
\item $0 \leq \weight(\beta)-\weight(\alpha) \leq 1$.
\end{enumerate}
\end{rem}

Any element $\alpha$ of $\wp{n}$ can be described by $\alpha=\{(A_1,w_1),\dotsc,(A_m,w_m) \}$ where $\{A_1,\dotsc A_r\}$ is a partition of $\{1,\dotsc,n\}$ and $w_i=\weight(A_i)$. We observe that $\wp{n}$ is a pure poset, i.e.~all maximal chains have the same length. 

\begin{figure}
\[
\xymatrix@R=50pt@C=15pt{ 
&*+[F-,]{123^0} \ar@{-}[dl]\ar@{-}[d] \ar@{-}[dr]  & &
*+[F-,]{123^1} \ar@{-}[dl] \ar@{-}[dll] \ar@{-}[dlll] 
                    \ar@{-}[dr] \ar@{-}[drr] \ar@{-}[drrr] & &
*+[F-,]{123^2} \ar@{-}[dl] \ar@{-}[d] \ar@{-}[dr] & \\
*+[F-,]{1^0|23^0} \ar@{-}[drrr] &
*+[F-,]{2^0|13^0} \ar@{-}[drr]  &
*+[F-,]{3^0|12^0} \ar@{-}[dr] & &
*+[F-,]{1^0|23^1} \ar@{-}[dl] &
*+[F-,]{2^0|13^1} \ar@{-}[dll] &
*+[F-,]{3^0|12^1} \ar@{-}[dlll] \\
& & & *+[F-,]{1^0|2^0|3^0} & && }
\]
\caption{\label{wpfigure}The poset $\wp{3}$}
\end{figure}

\begin{rem}
In Figure \ref{wpfigure}.~the weight $w$ of a block  $B=\{b_1,\dotsc,b_{k} \}$ is indicated by  $b_1\dotsb b_k^w$. E.g.~the block $\{1,2\}$ with weight $1$ is denoted by $12^1$.
\end{rem}

\begin{prop}
\label{wpisocomtwo}
The poset $\Pi_{\subtc{\Com}}(n)$ is isomorphic to $\wp{n}$.
\end{prop}
\begin{proof}
There is an obvious bijection between the elements of $\Pi_{\subtc{\Com}}(n)$ and $\wp{n}$ where a block $B$ enriched with an element $D^{|B|}_{i}$ with $i$ black product(s) corresponds to the same block $B$ with weight $i$ in $\wp{n}$. 

Now let $\alpha=\{(A_1,p_1),\dotsc,(A_m,p_m) \}$ be a $\Comtwo$-partition, then $\beta$ covers $\alpha$ iff 
\[
\beta=\{(A_j\cup A_k,\mu(\Bsmall;p_j,p_k)),(A_1,p_1),\dotsc,\widehat{(A_j,p_j)},\dotsc,\widehat{(A_k,p_k)},\dotsc,(A_m,p_m) \}
\]
or
\[
\beta=\{(A_j\cup A_k,\mu(\Wsmall;p_j,p_k)),(A_1,p_1),\dotsc,\widehat{(A_j,p_j)},\dotsc,\widehat{(A_k,p_k)},\dotsc,(A_m,p_m) \} .
\]
The first case corresponds to increasing the weight by one when merging two blocks of a weighted partition and the second case to keeping it constant, which precisely is the covering relation of $\wp{n}$.
\end{proof}

\subsection{Proof of Koszulness}
\label{koszulness} 

\begin{prop}
\label{comtwoCLlemma}
Let $\Pcal$ be a weakly associative binary quadratic set operad such that the maximal intervals of $\Pi_\Pcal$ are totally semimodular. Then the maximal intervals of $\Pi_{\subtc{\Pcal}}$ are CL-shellable.
\end{prop}

\begin{proof}
By Propositions \ref{posetoftc} and \ref{wpisocomtwo} we have that $\Pi_{\subtc{\Pcal}}=\Pi_\Pcal\times_\Pi \Pi_{\subtc{\Com}}=\Pi_\Pcal\times_\Pi\wp{{}}$. By Theorem 3.2 of \cite{Bjorner1983} CL-shellable is equivalent to admitting a recursive atom ordering. We aim to show that $\Pi_\Pcal\times_\Pi\wp{{}}$ admits such an ordering.

When denoting decorated partitions we will suppress the blocks only containing one element e.g.
\begin{align*}
\{(\{i,j\},p),(\{k,l\},p')\}=\{ & (\{i,j\},p),(\{k,l\},p'),(\{1\},1),\dotsc, \widehat{(\{i\},1)},\dotsc,\\
& \widehat{(\{j\},1)},\dotsc, \widehat{(\{k\},1)},\dotsc, \widehat{(\{l\},1)},\dotsc,(\{n\},1) \}.
\end{align*}

Denote the maximal elements $\{([n],p,w)\}$ of $\Pi_\Pcal(n)\times_{\Pi_n}\wp{n}$ by $\maximal{p}{w}$. Similarly denote the maximal elements $\{([n],p)\}$ of $\Pi_\Pcal(n)$ by $\maximalb{p}$. Assume that the length of a maximal interval $[\hat{0},\maximal{p}{w}]$ is greater than 1, otherwise we are done. We may also assume that the weight $w$ satisfies $0 < w < n-1$. Otherwise $[\hat{0},\maximal{p}{w}]$ is isomorphic to $[\hat{0},\maximalb{p}]\in\Pi_\Pcal(n)$ which is totally semimodular by assumption. Thus by Corollary 5.2 of \cite{Bjorner1983} it is CL-shellable.

Denote the atom $\{(\{i,j\},p)\}\in\Pi_\Pcal(n)$ by $\atomb{p}{i}{j}$. Similarly denote the atoms of $\Pi_\Pcal(n)\times_{\Pi_n}\wp{n}$ by $\atom{p}{w}{i}{j}$.

For a maximal interval $[\hat{0},\maximal{p}{w}]$ with $w>0$ we claim that any ordering of the form
\begin{align}
\label{atomorder}
&\atom{p_1}{0}{i_1}{j_1}\dashv \atom{p_1}{1}{i_1}{j_1}\dashv \dotsb \dashv
\atom{p_r}{0}{i_1}{j_1}\dashv \atom{p_r}{1}{i_1}{j_1}\dashv \dotsb \dashv 
\atom{p_1}{0}{i_m}{j_m}\dashv \atom{p_1}{1}{i_m}{j_m}\dashv \dotsb \dashv
\atom{p_r}{0}{i_m}{j_m}\dashv \atom{p_r}{1}{i_m}{j_m}
\end{align}
satisfies the second criterion \ref{recatomordering}\eqref{secondcriterion} of being a recursive atom ordering, where \newline 
$\{(\{i_k,j_k\},p_1),\dotsc,(\{i_k,j_k\},p_r)\}$ is some indexing of the atoms in $[\hat{0},\maximalb{p}]$ and $\alpha\dashv\beta$ means that $\alpha$ is less than $\beta$ in the atom ordering. Note that given an atom $\atom{p_s}{v}{i_t}{j_t}\in[\hat{0},\maximal{p}{w}]$, the atom $\atom{p_r}{\tilde{v}}{i_1}{j_1}$ will also be in the interval since we assume $0 < w < n-1$. Here $\tilde{v}$ denotes the element of $\{0,1\}\setminus\{v\}$.

Let $\atom{p_1}{w_1}{i}{j}$ and $\atom{p_2}{w_2}{k}{l}$ be distinct atoms with $\atom{p_1}{w_1}{i}{j}\dashv\atom{p_2}{w_2}{k}{l}$ and suppose $\atom{p_1}{w_1}{i}{j}$,$\atom{p_2}{w_2}{k}{l}\leq \gamma$, for some $\gamma=\{(C_1,q_1,v_1),\dotsc,(C_s,q_s,v_s)\}\in[\hat{0},\{([n],p,w)\}]$. We want to show that there is a $\delta\leq\gamma$ and an $\atom{p'}{w'}{i'}{j'}\dashv\atom{p_2}{w_2}{k}{l}$ such that $\atom{p'}{w'}{i'}{j'},\atom{p_2}{w_2}{k}{l}$ 
$\prec\delta$. Let $\gamma'=\{(C_1,q_1),\dotsc,(C_s,q_s)\}$. We have three main cases to consider.
\begin{enumerate}
\item $\{i,j\}$=$\{k,l\}$. Since the length of $[\hat{0},\{([n],p,w)\}]$ is greater than 1 we have that $n \geq 3$. We get two further subcases:

\begin{enumerate}
\item \label{subcaseA}$p_1=p_2$:
Since $\atom{p_1}{w_1}{i}{j}$,$\atom{p_1}{w_2}{i}{j}\leq \gamma$  and $w_1\neq w_2$ there must be at least one decorated block $(C_r,q,u)$ of $\gamma$ such that $\atomb{p_1}{i}{j}\leq\{(C_r,q)\}$ for some $q\in\Pcal_{|C_r|}$ and $|C_r|\geq 3$. Further, since $\Pcal$ is weakly associative there exist $q'\in\Pcal_3$ and $m\in C_r \setminus \{i,j\}$, for some $r$, such that $\delta'=\{(\{i,j,m\},q')\}\succ \atomb{p_1}{i}{j}$ and $\delta'\leq \gamma'$. Then $\delta=\{(\{i,j,m\},q',$ 
$\max(w_1,w_2))\}\leq\gamma$ and covers $\atom{p_1}{w_1}{i}{j}$ and $\atom{p_2}{w_2}{k}{l}$.

\item \label{subcaseB} $p_1\neq p_2$: 
Since $\Pi_\Pcal(n)$ is totally semimodular there exists a $\delta'=\{(\{i,j,m\},p)\}\in[\hat{0},\gamma']$ covering both $\atomb{p_1}{i}{j}$ and $\atomb{p_2}{k}{l}$. Then $\delta=\{(\{i,j,$ 
$m\},q,v)\}\leq\gamma$ covers both $\atom{p_1}{w_1}{i}{j}$ and $\atom{p_2}{w_2}{k}{l}$, where 
\[
v= \left\{ \begin{array}{ll}
               \max(w_1,w_2) & \mbox{if $w_1\neq w_2$}   \\
               w  & \mbox{if $w_1=w_2=w$}   \\
               w_1+1  & \mbox{if $w_1=w_2\leq w$}   \\
              \end{array}
       \right.
\]
\end{enumerate}

\item \label{caseB} $\{i,j\}\cap\{k,l\}=\{m\}$, for some $m\in\{i,j\}$. Let $m'$ be the element of $\{k,l\}\setminus\{m\}$. Since both atoms are less then $\gamma$ we must have that $\{i,j,m'\}$ is a subset of a block $C_r$ in $\gamma$. Since $\Pi_\Pcal(n)$ is totally semimodular there exists a $\delta'=\{(\{i,j,m'\},q)\}\in[\hat{0},\gamma']$ covering both $\atomb{p_1}{i}{j}$ and $\atomb{p_2}{k}{l}$.
Then $\delta=\{ (\{i,j,m'\},q,v) \}$ is an element covering both $\atom{p_1}{w_1}{i}{j}$ and $\atom{p_2}{w_2}{k}{l}$ and which is less than $\gamma$, where $v$
is as in case (\ref{subcaseB}).

\item $\{i,j\}\cap\{k,l\}=\emptyset$. Here we have two subcases:
\begin{enumerate}
\item $w_1\neq w_2$: $\delta=\{(\{i,j\},p_1,w_1),(\{k,l\},p_2,w_2)\}$ covers both $\atom{p_1}{w_1}{i}{j}$ and $\atom{p_2}{w_2}{k}{l}$ and will always be less than or equal to any $\gamma$ greater than both atoms. 
\item $w_1=w_2$: By the ordering of the atoms $\atom{p_1}{w_1}{i}{j}\dashv\atom{p_2}{w_2}{k}{l}$ implies $\atom{p_1}{\tilde{w_1}}{i}{j}\dashv\atom{p_2}{w_2}{k}{l}$, where $\tilde{w_1}$ is the element in $\{0,1\}\setminus \{w_1\}$. Now since  $\atom{p_1}{w_1}{i}{j}$,$\atom{p_2}{w_2}{k}{l}\leq \gamma$ either $\delta=\{(\{i,j\},p_1,w_1),(\{k,l\},p_2,w_2)\}\leq\gamma$ or $\tilde{\delta}=\{(\{i,j\},p_1,\tilde{w_1}),(\{k,l\},p_2,w_2)\}\leq\gamma$, where $\delta$ covers $\atom{p_1}{w_1}{i}{j}$ and $\atom{p_2}{w_2}{k}{l}$ whereas $\tilde{\delta}$ covers $\atom{p_1}{\tilde{w_1}}{i}{j}$ and $\atom{p_2}{w_2}{k}{l}$.
\end{enumerate}

\end{enumerate}

We also need to show that, given an ordering of the form (\ref{atomorder}), any interval $[\atom{q}{v}{i}{j},\maximal{p}{w}]$ satisfies the first criterion \ref{recatomordering}\eqref{firstcriterion} of being a recursive atom ordering. Note that $\atomb{q}{i}{j}\leq\maximalb{p}$ implies that there exists a $p'\in\Pcal_{n-1}$ such that $p=\mu(p';q,1,\dotsc,1)$ whence we observe that $[\atom{q}{v}{i}{j},\maximal{p}{w}] \iso [\hat{0},\maximal{p'}{w-v}]\subset \Pi_\Pcal(n-1)\times_{\Pi_{n-1}}\wp{n-1}$.

Thus checking the above step is readily done if we may order the atoms of $[\atom{q}{v}{i}{j},\maximal{p}{w}]$ in the same way as above. We only need to show that some way of ordering the atoms in pairs as above satisfies that the first atoms are the ones covering some atom $\atom{p'}{w'}{i'}{j'}\dashv\atom{p}{w}{i}{j}$. After that we can proceed by induction. 

We may assume that the length of $[\atom{q}{v}{i}{j},\maximal{p}{w}]$ is greater than 1, since otherwise we are done. We may also assume that $0 < w-v < n-2$, since otherwise the interval $[\atom{q}{v}{i}{j},\maximal{p}{w}]$ is isomorphic to $[\atomb{q}{i}{j},\maximalb{p}]\subset\Pcal_n$ which is totally semimodular by assumption.  Thus in the case that $w-v=0$ or $w-v=n-2$ we have by Theorem 5.1 of \cite{Bjorner1983} that any ordering of the atoms is a recursive atom ordering. We would therefore be able to freely order the atoms of $[\atom{q}{v}{i}{j},\maximal{p}{w}]$ so that the atoms that come first are those that cover some atom less than $\atom{q}{v}{i}{j}$ in the ordering (\ref{atomorder}).

Now the atoms are either of the form $\{(\{i,j\},q,v),(\{k,l\},t,u)\}$ which we denote by $\beta^{t,u}_{k,l}$ or of the form $\{(\{i,j,k\},t,v+u) \}$ which we denote by $\beta^{t,u}_{k}$, where $u\in\{0,1\}$. Let $\tilde{u}$ be the element of $\{0,1\}\setminus\{u\}$.

We have that  $\beta^{t,u}_{k,l}$ covers some $\atom{q'}{v'}{i'}{j'}\dashv\atom{q}{v}{i}{j}$, namely $\atom{q'}{v'}{i'}{j'}=\atom{t}{u}{k}{l}$,  iff $\atom{t}{u}{k}{l} \dashv\atom{p}{w}{i}{j}$. Since by the atom ordering of $[\hat{0},\maximal{p}{w}]$ we have that $\atom{t}{u}{k}{l} \dashv \atom{q}{v}{i}{j}$ iff $\atom{t}{\tilde{u}}{k}{l} \dashv \atom{q}{v}{i}{j}$, we have that  $\beta^{t,u}_{k,l}$ covers some $\atom{q'}{v'}{i'}{j'} \dashv \atom{q}{v}{i}{j}$ iff $\beta^{t,\tilde{u}}_{k,l}$ covers some $\atom{q'}{v'}{i'}{j'}\dashv\atom{q}{v}{i}{j}$.

Similarly we have that $\beta^{t,u}_{k}$ may cover some $\atom{q'}{v'}{i'}{j'}\dashv\atom{q}{v}{i}{j}$, where $\{i',j'\}\subset\{i,j,k\}$ and $q'$ is an appropriate element of $\Pcal_2$. Again $\atom{q'}{v'}{i'}{j'} \dashv \atom{p}{w}{i}{j}$ iff $\atom{q'}{\tilde{v'}}{i'}{j'} \dashv \atom{q}{v}{i}{j}$. Hence $\beta^{t,u}_{k}$ covers some  $\atom{q'}{v'}{i'}{j'}\dashv\atom{q}{v}{i}{j}$ iff $\beta^{t,\tilde{u}}_{k}$ does.

Thus we may order the atoms of $[\atom{q}{v}{i}{j},\maximal{p}{w}]$ by first putting all pairs of atoms, differing only in weight, covering some atom less than $\atom{q}{v}{i}{j}$ followed by all pairs of atoms not covering any atom less than $\atom{q}{v}{i}{j}$. Using the aforementioned identification $[\atom{q}{u}{i}{j},\maximal{p}{w}] \iso [\hat{0},\maximal{p'}{w-u}]$, we proceed by induction.

\end{proof}

\begin{rem}
Note that both the assumption that $\Pcal$ is weakly associative as well as the assumption that the maximal intervals of its associated poset are totally semimodular are necessary for the proof to go through. Both are used, in subcase \eqref{subcaseA} and e.g.~subcase \eqref{subcaseB}, respectively, and neither of the two properties implies the other.
\end{rem}

\begin{thm}
Let $\Pcal$ be a weakly associative binary quadratic basic-set operad such that the maximal intervals of $\Pi_\Pcal$ are totally semimodular, then $\tc{\aao{\Pcal}}$ and $\lc{(\aao{\Pcal}^{!})}$ are Koszul.
\end{thm}
\begin{proof}
By Proposition \ref{comtwoCLlemma} and the chain of implications \eqref{chainofimplications} we obtain that the maximal intervals of $\Pi_{\subtc{\Pcal}}$ are Cohen-Macaulay and by Corollary \ref{tcbasicset} we see that $\tc{\Pcal}$ is basic-set. Thus we can apply Theorem \ref{Bruno-main-thm} and conclude that $\tc{\aao{\Pcal}}$ is Koszul. Since $\tc{\aao{\Pcal}}$ and $\lc{(\aao{\Pcal}^{!})}$ are Koszul dual to each other we are done.
\end{proof}

We get the following immediate corollary.

\begin{cor}The following operads are Koszul: $\Comtwo$, $\Lietwo$, $\Permtwo$, $\Prelietwo$, $\Comtriastwo$, $\Postlietwo$, $\Astwotc$, $\Astwolc$, $\Diastwo$, $\Dendtwo$, $\Triastwo$ and $\Tridendtwo$.
\end{cor}
\begin{proof}
The operads $\Com$, $\Perm$, $\Comtrias$, $\As$, $\Dias$, and $\Trias$ are all algebraic operads which are linearizations of weakly associative basic-set operads whose associated posets have totally semimodular maximal intervals. The other operads are their Koszul dual operads. See \cite{Vallette2007,Chapoton2006} for these results and definitions of the operads.
\end{proof}

\section*{Acknowledgements}
The author is grateful to B.~Vallette for suggesting the use of fiber, Hadamard, black and white products in the description of compatible structures as well as for useful comments on the manuscript. The author wishes to thank also S.~Merkulov for valuable comments on the manuscript.

This note was typeset using  \Xy-pic.

\bibliographystyle{alpha}
\bibliography{compatiblestructures}

\end{document}